

\baselineskip=14pt
\parskip=10pt
\def\halmos{\hbox{\vrule height0.15cm width0.01cm\vbox{\hrule height
  0.01cm width0.2cm \vskip0.15cm \hrule height 0.01cm width0.2cm}\vrule
  height0.15cm width 0.01cm}}

\magnification=\magstephalf

\def\1{{\overline{1}}}
\def\2{{\overline{2}}}
\parindent=0pt
\overfullrule=0in

\def\frac#1#2{{#1 \over #2}}
\centerline
{\bf 
A quick proof that  $321$-avoiding permutations without double deficiencies 
}
\centerline
{\bf
are counted by the Motzkin numbers
}

\bigskip
\centerline
{\it  Tipaluck KRITYAKIERNE, Thotsaporn ``Aek'' THANATIPANONDA, and Doron ZEILBERGER}

$\bullet$ A permutation $\pi$ of length $n$ is  {\it $321$-avoiding}  if there are no triples $1 \leq i <j<k \leq n$ such that $\pi(i) >\pi(j)>\pi(k)$. They are famously counted by the
Catalan numbers, OEIS sequence $A108$.

$\bullet$ $k$, $1 \leq k \leq n$,  is a {\it double deficiency} of a permutation $\pi$ (of length $n$),  if $\pi(k)<k<\pi^{-1}(k)$.

$\bullet$  $w_1 \dots w_n \in \{ -1,0,1\}^n$  is a {\it Motzkin word}, if $\sum_{i=1}^{j} w_i$ is $\geq 0$ if $j<n$ and $0$ if $j=n$.
The number of such words is the $n$-th Motzkin number, written $M_n$ (OEIS sequence $A1006$).

The theorem below was first proved in [1] by a bijection whose statement was rather complicated and opaque. Here we present a more transparent proof. Note that the bijection is
the same, only its description is far simpler.

{\bf Theorem ([1]):} The number of $321$-avoiding permutations of length $n$ without double deficiencies is $M_n$

{\bf Proof}:  {\bf I. Motzkin to NoDD321}: Given a Motzkin word $w$, define increasing lists of integers, $I$ and $V$, as follows:

$i \in I$ if and only if $w_i = 0$ or $1$;  $i \in  V$ if and only if $w_i = 0$ or $-1$.

Let $\overline{I}:=\{1, \dots n \} \backslash I$ and $\overline{V}:=\{1, \dots, n\} \backslash V$  be the sorted complementary lists.
To get the output  $\pi$,
we define $ \pi(I(s)):= V(s) \, , \,  s=1,..,|I|$, and  $ \pi(\overline{I}(s)):= \overline{V}(s) \, , \,  s=1,..,n-|I|$.

Clearly $\pi$ is 
$321$-avoiding permutations where $V$ is the set of left-to-right maxima.
Since $I \, \cup V = \{1,2,\dots,n\}$, $\pi$ has no double deficiencies.

{\bf Example:} Consider the Motzkin word $w=[0,1,1,-1,0,1,-1,0,-1,1,-1]$.
Then $I = [ 1,2,3,5,6,8,10]$,  $V =[1,4,5,7,8,9,11]$ ; $\overline{I} = [4,7,9,11]$,  $\overline{V} = [2,3,6,10]$ , getting:
$\pi = [1,4,5,2,7,8,3,9,6,11,10] .$

{\bf  II. NoDD321 to Motzkin}: Let $\pi$ be a $321$-avoiding permutation without double deficiencies.
Let $V$ be the set of left-to-right maxima,
and let $I$ be the set of  $i$ such that $\pi(i) = v$ for some $v \in V$.
The pair $(I,V)$ determines a unique Motzkin word $w$ as follows:
$$
w_i =
\cases{
0, & if  $i \in I \cap V$ ; \cr
1, & if  $i \in I \cap \overline{V}$   ; \cr
-1, & if  $ i \in \overline{I} \cap V$ . \cr
}
$$
Since $\pi$ has no double-deficiencies, 
$I \, \, \cup \, V = \{1,2,\dots,n\}$.
By the condition of left-to-right maxima,
the resulting word has all its partial sums non-negative. $\halmos$
\vfill\eject

{\bf Reference}

[1] Martin Rubey and Christian Stump, {\it Double Deficiencies of Dyck Paths via the Billey-Jockusch-Stanley Bijection}, Journal of Integer Sequences, Vol. 20 (2017), Article 17.9.6 \hfill\break
{\tt https://cs.uwaterloo.ca/journals/JIS/VOL20/Stump/stump7.html}

\bigskip
\hrule
\bigskip
Tipaluck Krityakierne, Department of Mathematics, Faculty of Science, Mahidol University, Bangkok 10400, Thailand \hfil\break
Email: {\tt  tipaluck.kri at mahidol dot edu}
\bigskip
Thotsaporn ``Aek'' Thanatipanonda, Science Division, Mahidol University International College, Nakhon Pathom 73170, Thailand \hfil\break
Email: {\tt thotsaporn at gmail dot com}
\bigskip
Doron Zeilberger, Department of Mathematics, Rutgers University (New Brunswick), Hill Center-Busch Campus, 110 Frelinghuysen
Rd., Piscataway, NJ 08854-8019, USA. \hfill\break
Email: {\tt DoronZeil at gmail  dot com}   \quad .
\bigskip

{\bf July 2,  2026} \quad .

\end